\documentclass[11pt,a4paper]{article}
\usepackage{amsmath,amscd,amssymb,amsfonts,amsthm}
\usepackage[mathscr]{eucal}
\usepackage{hyperref}

\newtheorem{thm}{Theorem}[section]

\newtheorem{prop}{Proposition}[section]

\newcommand{\Sch}{{\mathfrak{S}}}

\newcommand{\p}{\partial}

\def\a{\alpha}

\def\D{\Delta}

\renewcommand{\S}{{\Sigma}}

\def\t{\theta}

\newcommand{\bs}{{\boldsymbol{s}}}

\numberwithin{equation}{section}

\usepackage{amsmath,amscd,amssymb,amsfonts}
\usepackage[mathscr]{eucal}
\usepackage{hyperref}

\newtheorem{df}{Definition}[section]

\theoremstyle{remark}

\numberwithin{equation}{section}

\renewcommand{\a}{{\alpha}}
\renewcommand{\epsilon}{{\varepsilon}}
\newcommand{\brho}{{\boldsymbol \rho}}
\newcommand{\bsigma}{{\boldsymbol \sigma}}

\renewcommand{\S}{{\mathcal S}}
\newcommand{\PTM}{{\Pi T^*M}}
\renewcommand{\t}{\tilde}

\begin{document}

\centerline {\bf Laplacians in  Odd Symplectic Geometry}

\bigskip

\centerline {\bf Hovhannes M.~Khudaverdian}

\bigskip

\centerline {\it Department of Mathematics, UMIST,
 Manchester M60 1QD, UK}

\centerline {\it and LIT, JINR, Dubna, 141980, Russia.}

\centerline {\it e-mail: khudian@umist.ac.uk}

%\subjclass[2000]{Primary 53D05, 58A50. Secondary 53D17, 58J60,81T70.}

\begin{abstract}
We consider odd Laplace operators arising in odd symplectic
geometry. Approach based on semidensities (densities of weight
$1/2$) is developed. The role of semidensities in the
Batalin--Vilkovisky formalism is explained. In particular, we
study the relations between semidensities on an odd symplectic
supermanifold and differential forms on a purely even Lagrangian
submanifold. We establish a criterion of ``normality'' of a volume
form on an odd symplectic supermanifold in terms of the canonical
odd Laplacian acting on semidensities.
\end{abstract}

\tableofcontents

\section {Symplectic and Poisson structures}

 A symplectic structure on a manifold $M$ is defined by a
 non-degenerate closed two-form $\omega$. %$(M^{2n},\omega)$
 In a vicinity of an arbitrary point one can consider
 coordinates $(x^1,\dots,x^{2n})$ such that $\omega=\sum_{i=1}^n dx^i dx^{i+n}$.
 Such coordinates are called Darboux coordinates.
 To a symplectic structure corresponds a non-degenerate Poisson structure $\{\,\,,\,\,\}$.
 In Darboux coordinates $\{x^i,x^j\}=0$ if $|i-j|\not=n$ and
 $\{x^i,x^{i+n}\}=-\{x^{i+n},x^{i}\}=1$.
 The condition of closedness of the two-form $\omega$ corresponds to the Jacobi identity
 $\{f,\{g,h\}\}+\{g,\{h,f\}\}+\{h,\{f,g\}\}=0$ for the Poisson bracket.
 If a symplectic or Poisson structure is given, then
 every function $f$ defines a vector field (the Hamiltonian vector field)
 ${\bf D}_f$ such that ${\bf D}_f g=\{f,g\}=-\omega ({\bf D}_f,{\bf D}_g)$.

A Poisson structure can be defined   independently of a symplectic
structure (see below).
 In general it can be degenerate, i.e., there exist non-constant functions
 $f$ such that ${\bf D}_f=0$. In the case when a Poisson structure
 is non-degenerate (corresponds to a symplectic structure),
 the map from $T^*M$ to $TM$ defined by the relation
   $f\mapsto {\bf D}_f$ is an isomorphism.

One can straightforwardly generalize these constructions to the
    supercase and
    consider  symplectic and Poisson structures
    (even or odd) on supermanifolds.
    An even (odd) symplectic structure on a supermanifold is defined
    by an even (odd) non-degenerate closed two-form.
    In the same way as the  existence of a symplectic structure
    on an ordinary manifold implies that the manifold is even-dimensional
    (by the non-degeneracy condition for the form $\omega$),
    the existence of an even or odd symplectic structure
    on  a supermanifold implies that the dimension of the supermanifold
    is equal either to $(2p.q)$ for an even structure or to  $(m.m)$ for
    an odd structure. Darboux coordinates exist in both cases.
    For an even structure, the two-form in Darboux coordinates
    $z^A=(x^1,\dots,x^{2p};\theta_1,\dots,\theta_q)$
    has the form $\sum_{i=1}^p dx^idx^{p+i}$
    $+\sum_{a=1}^q \varepsilon_a d\theta_a d\theta_a$, where
    $\varepsilon_a=\pm 1$.
   For an odd structure, the two-form in Darboux coordinates
    $z^A=(x^1,\dots,x^{m};\theta_1,\dots,\theta_m)$
    has the form $\sum_{i=1}^m dx^id\theta_{i}$.

   The  non-degenerate odd Poisson bracket corresponding
   to an odd symplectic structure has the following
   appearance in Darboux coordinates: $\{x^i, x^j\}=0$, $\{\theta_i,\theta_j\}=0$
   for all $i,j$ and
   $\{x^i,\theta_j\}=-\{\theta_j,x^i\}=\delta^i_j$. Thus for
   arbitrary two
   functions $f,g$
                  \begin{equation}\label{oddpoissondefinition}
     \{f,g\}=\sum_{i=1}^n\left(
      \frac{\p f}{\p x^i}
     \frac{\p g}{\p \theta_i}+
                 (-1)^{p(f)}
     \frac{\p f}{\p \theta_i}
      \frac{\p g}{\p x^i}
     \right),
                 \end{equation}
where we denote by $p(f)$ the parity of a function $f$ (e.g.,
$p(x^i)=0$, $p(\theta_j)=1$). Similarly one can write down the
formulae for the non-degenerate even Poisson structure
corresponding to an even symplectic structure.

A Poisson structure (odd or even) can be defined on a
supermanifold  independently  of a symplectic structure as a
bilinear operation on functions (bracket) satisfying the following
relations taken as axioms:
%\begin{equation}
\begin{gather}
p\left(\{ f,g\}\right)=p(f)+p(g)+\epsilon,\\
    \{f,g\}=-\{g,f\}(-1)^{(p(f)+\epsilon)(p(g)+\epsilon)}, \\
          \{f,gh\}=\{f,g\}h+\{f,h\}g(-1)^{p(g)p(h)} \text{\quad(Leibniz rule)},\\
    \{f,\{g,h\}\}(-1)^{(p(f)+\epsilon)(p(h)+\epsilon)}+
            \text{cycl.} =0
             \text{\quad(Jacobi identity)},
\end{gather}
%\end{equation}
where $\epsilon$ is the parity of the bracket ($\epsilon=0$ for an
even Poisson structure and
  $\epsilon=1$ for an odd one).
  The correspondence between functions and Hamiltonian vector
  fields is defined in the same way  as on ordinary manifolds: ${\bf D}_f g=\{f,g\}$. Notice a
possible parity shift: $p({\bf D}_f)=p(f)+\
   \epsilon$.
   Every Hamiltonian vector field
  ${\bf D}_f$ defines an infinitesimal transformation
  preserving the Poisson structure (and the corresponding symplectic structure
  in the case of a non-degenerate Poisson bracket).

   Notice that even or odd Poisson structures on an arbitrary supermanifold
    can be obtained as ``derived'' brackets from the canonical
   symplectic structure on the cotangent bundle, in the following way.

  Let $M$ be a supermanifold and  $T^*M$ be its cotangent bundle. By changing parity
  of coordinates in the fibres of $T^*M$ we arrive at the supermanifold
   $\Pi T^*M$. If $z^A$ are arbitrary coordinates on the supermanifold $M$, then
   we denote by $(z^A,p_{B})$ the corresponding coordinates on the supermanifold
   $T^*M$ and by $(z^A, z_{*B})$ the corresponding coordinates on $\Pi T^*M$:
  $p(z^A)=p(p_{A})=p(z_{*A})+1$. If $(z^{A'})$ are another coordinates
  on $M$, $z^A=z^A(z')$,  then the coordinates $z_{*A}$ transform
  in the same way as the coordinates $p_{A}$ (and as the partial derivatives $\p/\p z^A$):
                              \begin{equation}
                        \label{transformationincotangentbundle}
             p_{A'}= \frac{\p z^B(z')}{\p z^{A'}}\,p_{B}
              \text{\quad and \quad}
              z_{*A'}=\frac{\p z^B(z')}{\p z^{A'}}\,z_{*B}.
                               \end{equation}
    One can consider the canonical non-degenerate even
   Poisson structure $\{\,\,,\,\,\}_0$ (the canonical even symplectic structure) on $T^*M$
    defined by the relations
    $\{z^A,z^B\}_0=\{p_{C},p_{D}\}_0=0$, $\{z^A,p_{B}\}_0=\delta^A_B$,
   and, respectively, the canonical non-degenerate odd
   Poisson structure $\{\,\,,\,\,\}_1$ (the canonical odd symplectic structure) on $\Pi T^*M$
   defined by the relations
    $\{z^A,z^B\}_0=\{z_{*C},z_{*D}\}_0=0$, $\{z^A,z_{*B}\}_0=\delta^A_B$.

   Now consider Hamiltonians on $T^*M$ or on $\Pi T^*M$ that are quadratic
   in coordinates of the fibres. An arbitrary odd quadratic  Hamiltonian  on $T^*M$
     (an arbitrary even quadratic Hamiltonian on $\Pi T^*M$):
\begin{equation}%\label{}
     \Sch(z,p)=\Sch^{AB}\,p_{A}p_{B}  \ \ (p(\Sch)=1) \text{\quad or \quad}
     \Sch(z,z_*)=\Sch^{AB}z_{*A}z_{*B} \ \ (p(\Sch)=0),
\end{equation}
satisfying the
     condition that the canonical Poisson bracket of this Hamiltonian
     with itself vanishes:
\begin{equation}\label{masterodd}
\{\Sch,\Sch\}_0 =0 \text{\quad or\quad} \{\Sch,\Sch\}_1=0\,
\end{equation}
defines an odd Poisson structure (an even Poisson structure) on
$M$ by the formula
                         \begin{equation}\label{bracketfrommasterhamiltonian}
     \{f,g\}_{\epsilon+1}^\Sch=\{f,\{\Sch,g\}_{\epsilon}\}_{\epsilon}\,.
                         \end{equation}
   The Hamiltonian $\Sch$ which generates an odd (even) Poisson structure on $M$
   via the canonical even (odd) Poisson structure on $T^*M$ ($\Pi T^*M$)
   can be called the master Hamiltonian. The bracket \eqref{bracketfrommasterhamiltonian}
   is a ``derived bracket''.
 The Jacobi identity for it is equivalent to the vanishing of the canonical Poisson bracket
 for the master Hamiltonian.
 One can see that an arbitrary Poisson structure on a supermanifold
 can be obtained as a derived bracket.

 What happens if we change the parity of the master Hamiltonian in
 \eqref{bracketfrommasterhamiltonian}?
 The answer is the following. If $\Sch$ is an even quadratic Hamiltonian
 on $T^*M$ (an odd quadratic Hamiltonian on $\Pi T^*M$), then the
 condition of vanishing of the canonical even Poisson bracket $\{\,\,,\,\,\}_0$
 (the canonical odd Poisson bracket $\{\,\,,\,\,\}_1$)  becomes
 empty (it is obeyed automatically) and the relation \eqref{bracketfrommasterhamiltonian}
 defines an even Riemannian metric (an odd Riemannian metric) on $M$.

Formally, odd symplectic (and odd Poisson) geometry  is a
generalization of symplectic (Poisson) geometry to the supercase.
However, there are unexpected analogies between the constructions
in odd symplectic geometry and in Riemannian geometry (see
\cite{sdens:voron} and later below).  The construction of derived
brackets could explain close relations between odd Poisson
structures in supermathematics and the Riemannian geometry (see
\cite{sdens:voron}).

The construction of the derived bracket (without the name) and the
elaboration of the unified viewpoint for different geometries in
terms of derived brackets are due to T.~Voronov~\cite{tv:private}.
Derived brackets (under this name) were independently introduced
and studied in~\cite{yvette:derived}. It has to be noted that in
the physical literature the relations of the type
\eqref{bracketfrommasterhamiltonian} for brackets of different
parity were considered in \cite{ners0:brackets} and
\cite{ners:brackets}, where they were used for obtaining derived
brackets on Lagrangian surfaces.  This approach was considered
later in \cite{maks:derivedbrackets} and essentially developed in
\cite{bf:derivedbrackets}, where  constructions involving
generalized ``higher order'' even and odd Poisson brackets
appeared.

In what follows we consider second order differential operators on
an odd symplectic supermanifold and study their geometric
properties.   Some of our constructions can be automatically
considered in the  case of a general odd Poisson structure.

\section {Odd Laplacians on functions}

\subsection {Definition and properties}\label{subsecdefinitionandproperties}

In ordinary symplectic geometry the symplectomorphisms of $M^{2n}$
(the diffeomorphisms preserving the two-form $\omega$) preserve
the volume form $\brho_\omega=\omega^n$ (Liouville's theorem).
What is the situation in the supercase?  In spite of the fact that
differential forms are not objects of integration on
supermanifolds (see details in~\cite{tv:git}), the Liouville
theorem still holds in even symplectic geometry. One can see that
the coordinate volume form $\brho={\mathcal D}z$ in Darboux
coordinates defines a global volume form that is preserved under
symplectomorphisms. The situation is drastically different for an
odd symplectic structure. Let
$z^A=(x^1,\dots,x^n,\theta_1,\dots,\theta_n)$ be Darboux
coordinates. Consider, for example,  the  transformation
$x^1\mapsto 2x^1$, $\theta_1\mapsto \frac{1}{2}\,\theta_1$ to
another Darboux coordinates. The Berezinian of this transformation
is equal to $4$, hence the coordinate volume ${\mathcal D}z$ form
is not preserved. One can prove (see below) that on an odd
symplectic supermanifold  there is no volume form invariant w.r.t.
all symplectomorphisms.

Let $\brho$ be an arbitrary  volume form on an odd symplectic
supermanifold. Consider the linear differential operator
$\Delta_\brho$ on functions such that its action on a function $f$
is equal (up to a coefficient) to the divergence ${\rm div}_\brho
{\bf D}_f$ of the Hamiltonian vector field ${\bf D}_f$ w.r.t. to
the volume form $\brho$:

\begin{equation}\label{defoflaplacian}
                           \Delta_\brho f:= \frac{1}{2}\,(-1)^{p(f)}
                          {\rm div}_\brho {\bf D}_f=
                 \frac{1}{2}\,(-1)^{p(f)}{\mathcal L}_{{\bf D}_f}\log\brho=
                 \frac{1}{2}\,(-1)^{p(f)}\frac{{\mathcal L}_{{\bf D}_f}\brho}{\brho}\,.
\end{equation}

Let $(x^1,\dots,x^n;\theta_1,\dots,\theta_n)$ be Darboux
coordinates  and let $\brho=\rho(x,\theta){\mathcal D}(x,\theta)$
in these coordinates. It follows from \eqref{oddpoissondefinition}
that
\begin{equation*}%\label{}
    {\bf D}_f=\frac{\p f}{\p x^i}{\p\over\p\theta_i}+
                (-1)^{p(f)}\frac{\p f}{\p \theta_i}\frac{\p}{\p x^i}.
\end{equation*}
Hence
\begin{equation}\label{defoflaplaciana}
                        \Delta_\brho f=\Delta_0 f+\frac{1}{2}\,\{\log \rho, f\},
\end{equation}
where
\begin{equation}\label{defoflaplaciansimple}
                      \Delta_0 f=
        \frac{\p^2 f}{\p x^i\p\theta_i}.
\end{equation}
We come to the \textit{odd Laplacian on functions}, a second order
differential operator depending on the volume form.

Notice that \eqref{defoflaplacian} defines the operator
$\Delta_{\brho}$ in terms of the Poisson bracket. This expression
defines a linear operator on functions for an arbitrary Poisson
structure, even or odd. It is an operator of the second order for
an odd Poisson structure and an operator of the first order for an
even Poisson structure. In the even case this operator of the
first order is a Poisson vector field (the divergence of the
Poisson bivector) specifying the so-called Weinstein's modular
class of an even Poisson manifold~\cite{weinstein:modular}.

If $\brho$ and $\brho^\prime$ are two volume forms,
$\brho^\prime=g\brho$, then
\begin{equation}\label{changingoflaplacian}
                 {\rm div}_{\brho^\prime} f- {\rm div}_\brho f=
                 \frac{{\mathcal L}_{{\bf D}_f}g}{g}=\{f, \log g\}.
\end{equation}

In particular,  the existence of a ``canonical'' volume form
$\brho_\omega$ preserved under all symplectomorphisms  would imply
that the operator~\eqref{defoflaplacian} is a first order
differential operator given by the r.h.s.
of~\eqref{changingoflaplacian}, because~\eqref{defoflaplacian}
would evidently vanish for $\brho=\brho_\omega$. On the other
hand, by \eqref{defoflaplaciana} and \eqref{defoflaplaciansimple},
$\Delta_\rho$ is a second order differential operator. Thus we
have proved that \textit{on an odd symplectic supermanifold there
is no canonical volume form}.

This should be compared with the even Poisson situation where
Weinstein's modular class (see above) is the obstruction to the
existence of a volume form invariant under all Hamiltonian flows,
and it vanishes in the (even) symplectic case. We see that the
situation with the odd bracket (symplectic or not) is more
complicated.

Now consider the properties of the odd Laplacian $\Delta_\brho$.

One can see that for an arbitrary odd Poisson supermanifold   the
Leibniz rule for the second derivatives takes the following   form
for the odd Laplacian:
\begin{equation}\label{yieldingrelation}
          \Delta_\brho (f\cdot g)=\Delta_\brho f\cdot g+(-1)^{p(f)}\{f,g\}+
                    (-1)^{p(f)} f\cdot\Delta_\brho g\,.
\end{equation}
In other words, the operator $\Delta_\brho$ generates the Poisson
structure.

We have already mentioned that \eqref{yieldingrelation} has a
straightforward analogue in Riemannian geometry:   $\Delta
(fg)=g\Delta f +\langle df,dg\rangle +f \Delta g$, where $\langle
\,\,,\,\,\rangle $ is the scalar product given by the Riemannian
metric and $\Delta$ is the Beltrami--Laplace operator
corresponding to the metric (see details in~\cite{sdens:voron}).

Another very important property of the odd  Laplacian (on an
arbitrary odd Poisson supermanifold) is that it preserves the
Poisson bracket:
\begin{equation}\label{bracketpreservation}
             \Delta_\brho \{f, g\}=\{\Delta_\brho f, g\}+
                    (-1)^{p(f)+1} \{f,\Delta_\brho g\}.
\end{equation}

Now let us return to the  canonical odd symplectic  structure on
$\Pi T^*M$ (see the first section).   We consider the case when
$M$ is a usual $n$-dimensional manifold.  The base manifold $M$ is
a Lagrangian $(n.0)$-dimensional surface in the
$(n.n)$-dimensional odd symplectic supermanifold $\Pi T^*M$. This
example can be considered as the basic example of an odd
symplectic supermanifold \footnote{For a usual symplectic manifold
$E$ and a Lagrangian surface $L$ in it there exists a
symplectomorphism  between the cotangent bundle $T^*L$ and, in
general, only a tubular neighborhood  of $L$ in $E$. The
triviality of topology in odd directions allows to identify the
cotangent bundle $\Pi T^*L$ to a purely even Lagrangian surface in
an odd symplectic supermanifold $E$ with the whole supermanifold
$E$ if  $L$ coincides with the underlying manifold of $E$ (see
\cite{hov:semi} for details).}.

Functions on $\Pi T^*M$ encode multivector   fields on $M$:
\begin{multline}\label{identificationofvectors}
        f(x,\theta)=f(x)+f^i(x)\theta_i+f^{ik}(x)\theta_i\theta_k+\dots \\ \leftrightarrow
                T= f(x)+f^i(x)\p_i+f^{ik}(x)\p_i\wedge\p_k+\dots
\end{multline}
The odd Poisson bracket of functions on $\Pi T^*M$ corresponds to
the  Schouten  bracket  (``skew-symmetric concomitant'') of
multivector fields.

To every diffeomorphism of $M$ naturally corresponds the induced
symplectomorphism  of $\Pi T^*M$, but in general one can consider
symplectomorphisms that do not correspond to diffeomorphisms of
the base and destroy the cotangent bundle structure (see the next
section).

Let us  now analyze the meaning of an odd
Laplacian~\eqref{defoflaplacian} on $\Pi T^*M$.   Let $(x^i)$ be
arbitrary coordinates on $M$ and $(x^i,\theta_j)$  the
corresponding Darboux coordinates on $\Pi T^*M$. Let
$\bsigma={\mathcal D}(x)$ and $\brho={\mathcal D}(x,\theta)$ be
the coordinate volume forms on $M$ and $\Pi T^*M$ respectively.
Then the odd Laplacian $\Delta_\brho$ on $\Pi T^*M$ is given by
the formula \eqref{defoflaplaciansimple}. It is obvious from
\eqref{defoflaplaciansimple} that in this case the action of the
operator $\Delta_\brho$ on functions on $\Pi T^*M$ corresponds to
the divergence of multivector fields on $M$ with respect to the
coordinate volume form $\bsigma$. Every volume form on $M$ has a
local appearance as a coordinate volume form (in some local
coordinate system). On the other hand, it follows from
\eqref{transformationincotangentbundle} that the determinant of an
arbitrary coordinate transformation $x\mapsto x'= x'(x)$ on $M$ is
equal to the  square root of the Berezinian of the corresponding
coordinate transformation $(x,\theta)\mapsto (x',\theta')$:
\begin{equation}\label{squarerelations}
{\rm Ber}\,{\p( x',\theta')\over\p( x,\theta)}= {\rm Ber}
\begin{pmatrix}
  \frac{\p x^{i'}}{\p x^i} & *\\
  0 & \frac{\p x^i}{\p x^{i'}}
\end{pmatrix}
= \Biggl(\det \Bigl(\frac{\p  x^{i'}}{\p x^i}\Bigr)\Biggr)^2\,.
\end{equation}

Hence we come to important conclusions: (a) to every volume form
$\bsigma$ on $M$ corresponds a volume form $\brho=\bsigma^2$ on
$\Pi T^*M$; (b) ${\rm div}_\bsigma{T}=\Delta_{\bsigma^2}f$, where
we identify multivector fields on $M$ with functions on $\PTM$
by~\eqref{identificationofvectors}; (c) for  a volume form
$\brho=\bsigma^2$ on $\Pi T^*M$ holds  the condition
\begin{equation}\label{nilpotencycondition}
                     \Delta_{\brho}^2=0,
\end{equation}
because the square of the divergence operator on multivector
fields  equals  zero.

These relations make a bridge between odd symplectic  geometry and
classical vector calculus. They are closely related with the
geometric  meaning of the Batalin--Vilkovisky formalism (see
\cite{hov:khn2}).

%\smallskip
What can we say about  $\Delta_\brho^2$  in the general case?

For an arbitrary odd Poisson supermanifold the operator
$\Delta_\brho^2$ is a Poisson vector field: $\Delta_\brho^2
(fg)=(\Delta_\brho^2 f)g+f(\Delta_\brho^2 g)$ and $\Delta_\brho^2$
preserves the Poisson bracket. This follows from relations
\eqref{yieldingrelation}, \eqref{bracketpreservation}. Under a
change of volume form $\brho\mapsto\brho^\prime=g\brho$ this
Poisson vector field changes by a Hamiltonian vector field:
\begin{equation}\label{cocyclecondition}
         \Delta^2_{\brho^\prime}-\Delta^2_\brho={\bf D}_{H(\brho^\prime,\brho)}
         \text{\quad where\quad }
         H(\brho^\prime,\brho)=\frac{1}{\sqrt g}\,\Delta_\brho \sqrt g\,.
\end{equation}
This relation leads to a non-trivial groupoid
structure~\cite{sdens:voron}.

For an odd symplectic supermanifold one can always pick a volume
form $\brho$ such that $\Delta_\brho^2=0$, namely, as we have
shown above, one can identify the symplectic supermanifold with
$\Pi T^*L$ for an purely even Lagrangian surface $L$ by a suitable
symplectomorphism and then choose $\brho=\bsigma^2$, yielding
\eqref{nilpotencycondition}.  Hence, it follows from
\eqref{cocyclecondition}  that for odd Laplacians on a symplectic
supermanifold the operator $\Delta_\brho^2$ is always {\rm a
Hamiltonian vector field}.

\begin{df}
A volume form $\brho$ on an odd symplectic supermanifold $E$ is
called \textit{normal}  if in a vicinity of an arbitrary point
there exist Darboux coordinates $(x,\theta)$ such that $\brho$ is
the coordinate volume form in these Darboux coordinates:
$\brho={\mathcal D}(x,\theta)$.
\end{df}

The volume form $\brho=\bsigma^2$ on $\Pi T^*M$
in~\eqref{nilpotencycondition} is a normal volume form by
definition. If a volume form $\brho$ is normal, then
$\Delta_\brho^2=0$, by the definition of the odd
Laplacian~\eqref{defoflaplaciana}.  Does the Hamiltonian field
$\Delta_\brho^2$ has to be equal to zero for every volume form
$\brho$? As it follows from~\eqref{cocyclecondition}, the answer
is, generally, no.  Does the condition $\Delta_\brho^2=0$ implies
that $\brho$ is a normal volume form? We will give a detailed
analysis in the next sections, where we will study the canonical
odd Laplacian acting on semidensities.   Now two words about where
odd Laplacians come from.

\subsection {Where an odd Laplacian comes from}

Odd Laplacians have appeared in mathematical physics for the first
time around 1981 in the pioneer  works by Batalin and
Vilkovisky~\cite{bv:perv, bv:vtor, bv:closure} for the purpose of
constructing a Lagrangian version of the BRST quantization (the
``BV-formalism''). Batalin and Vilkovisky introduced an odd
Laplacian acting on functions on an odd symplectic supermanifold
as $\Delta_0=\frac{\p^2}{\p x^i\p\theta_i}$, where
$(x^i,\theta_i)$ are some Darboux coordinates. (The invariant
definition~\eqref{defoflaplacian} of an odd Laplacian
$\Delta_\rho$ depending on a volume form $\brho$ was given later
in \cite{hov:deltabest}.) The following very important properties
of this operator were fixed in their works
(see~\cite{bv:closure}). If $(x',\theta')$ are another Darboux
coordinates, and $\Delta_0$ and $\Delta_0^\prime$ denote the odd
Laplacians ~\eqref{defoflaplaciansimple} in the Darboux
coordinates $(x,\theta)$ and $(x^\prime,\theta^\prime)$
respectively, then
\begin{equation}\label{equivalenttoone}
          \Delta_0=\Delta_0'+
          \frac{1}{2}\,\{\log {\rm
          Ber}\,\frac{\p(x',\theta')}{\p(x,\theta)}, \ \}
\end{equation}
and
\begin{equation}\label{batalinindentity}
                  \Delta_0\left({\rm Ber}
            {\p(x',\theta')\over\p(x,\theta)}
\right)^{1/2}=0 \text{\quad (``Batalin--Vilkovisky identity'')}.
\end{equation}
The property \eqref{equivalenttoone} of the operator $\Delta_0$ is
closely related with the  invariant
expression~\eqref{defoflaplaciana} for $\Delta_\rho$.  The second
property~\eqref{batalinindentity},  which was stated in
\cite{bv:closure} (in a non-explicit form it also appeared
in~\cite{tyutin:unknown}), is highly non-trivial.   This identity
is deeply related with the existence of canonical odd Laplacian on
semidensities (see the next section).

The operator $\Delta_0$ was introduced in \cite{bv:perv} for the
purpose of formulating the so-called Batalin--Vilkovisky quantum
master equation, which is the equation
\begin{equation}\label{masterequation}
                        \Delta_0 \sqrt{f(x,\theta)}=
                        0
\end{equation}
on the function $f(x,\theta)=\exp {{i\S(x,\theta)\over \hbar}}$,
where $\brho={\mathcal D}(x,\theta)$ is a coordinate volume form
in the space of fields and antifields ($(x,\theta)$ are Darboux
coordinates). The master action $\S(x,\theta)$ defines a measure
element on Lagrangian surfaces corresponding to gauge choice. This
measure is gauge invariant  if the master equation is satisfied
(see \cite{bv:perv, bv:vtor, bv:closure}). Using the identity
$\Delta_0 \exp g= (\Delta_0 g+1/2\{g,g\})\exp g $,   we can
rewrite the quantum master equation~\eqref{masterequation} as
$-4\hbar \Delta_0 \S+\{\S,\S\}=0$ and taking $\hbar\to 0$ we
arrive at the Batalin--Vilkovisky classical master equation:
$\{\S,\S\}=0$. The geometrical meaning of the master equation was
studied in~\cite{hov:khn1, hov:khn2} and most notably by
A.~S.~Schwarz in~\cite{ass:bv}. In particular, the following
result was obtained. Suppose $\brho$ is a normal volume form. Then
there are implications:
\begin{equation}\label{nashiteoremy}
\text{the volume form $\brho^\prime=f(x,\theta)\brho$ is normal\ \
}\Rightarrow\ \ \Delta_\brho\sqrt f=0\ \ \Rightarrow\ \
\Delta_{\brho^\prime}^2=0\,.
\end{equation}
In this statement the master equation  is not formulated
invariantly, but it stands  between two invariant conditions. The
exact statement about the relations between the three conditions
in~\eqref{nashiteoremy} will be formulated in the next section in
the language of semidensities.

\section {Canonical odd Laplacian on semidensities}

\subsection {Definition of the canonical  Laplacian}

A density of weight $t$ on a supermanifold is a function of local
coordinates such that under a change of variables it is multiplied
by the $t$-th power of the Berezinian of transformation.      We
will  consider semidensities (densities of weight $t=\frac{1}{2}$)
on an odd symplectic supermanifold.

First of all, let us  consider again the supermanifold $\Pi T^*M$
for a usual manifold $M$, with the canonical symplectic structure.
Recall that  functions on $\Pi T^*M$ encode multivector fields on
$M$ (see~\eqref{identificationofvectors}). Our
claim~\cite{hov:semi} is that \textit{semidensities on $\Pi T^*M$
encode differential forms on $M$}.

Indeed, let $(x^1,\dots,x^n)$ be arbitrary local coordinates on
the manifold $M$. Let $(x^1,\dots,x^n$ $;
\theta_1,\dots,\theta_n)$ be the corresponding Darboux coordinates
on $\Pi T^*M$. In the same way as multivector fields on $M$ can be
identified with functions on $\Pi T^*M$, differential forms on $M$
can be identified  with functions on the supermanifold $\Pi TM$
obtained from the tangent bundle $TM$ by changing parity of
coordinates in the fibres. If $(x^1,\dots,x^n,\xi^1,\dots,\xi^n)$
are the coordinates on $\Pi TM$ corresponding to coordinates
$(x^1,\dots,x^n)$ on $M$ ($p(\xi^k)=1$), then a differential form
$\omega_0(x)+dx^i\omega_i(x)+\dots$ can be identified with the
function $\omega(x,\xi)=\omega_0(x)+\xi^i\omega_i(x)+\dots$ on
$\Pi TM$.

Recall  that given a volume form $\bsigma=\sigma(x){\mathcal
D}x=\sigma(x)dx^1\dots dx^n$ on $M$ the Hodge operator transforms
a multivector field  on $M$ corresponding to the function
$f(x,\theta)$ on $\Pi T^*M$ to the differential form on $M$
corresponding to the function $\omega(x,\xi)$ on $\Pi TM$, where
\begin{equation*}%\label{}
    \omega(x,\xi)=\int \exp(i\xi^i\theta_i)f(x,\theta)\,\sigma(x)\D\theta
\end{equation*}
(a ``Fourier transform''). (In a conventional language it is a
contraction of a top order form with a multivector, but the
language of integrals is more flexible.) It follows that without a
volume form, the Hodge operator acts on multivector densities of
weight $t=1$ on $M$ transforming them into differential forms on
$M$, and vice versa.

On the other hand, from~\eqref{squarerelations} it follows  that
under canonical transformations on $\Pi T^*M$ induced by changes
of coordinates of $M$ a density on $M$ transforms as a semidensity
on $\Pi T^*M$.

Hence we come to a 1-1-correspondence between differential forms
on $M$ (functions on $\Pi TM$) and semidensities on $\Pi T^*M$, as
follows\footnote {In these considerations we assume that the
manifold $M$ is orientable and an orientation is chosen.}:
\begin{equation}\label{map}
             \omega(x,\xi)\mapsto \bs=s(x,\theta)\sqrt{\D(x,\theta)}=
                              \left(
             \int \exp(-i\xi\theta)\omega(x,\xi)\D\xi
                             \right)\,
                              \sqrt{D(x,\theta)}
\end{equation}
(This map for the first time appeared in~\cite{ass:bv} in a
non-explicit way.) For example, if $M$ is a two-dimensional
manifold, then $f(x)\mapsto f(x)\theta_1\theta_2\,\sqrt
{\D(x,\theta)}$, $\omega_1dx^1+\omega_2dx^2\mapsto
(\omega_1\theta_2-\omega_2\theta_1)\,\sqrt {\D(x,\theta)}$,
$\omega dx^1dx^2$ $\mapsto -\omega\,\sqrt{\D(x,\theta)}$

The relation~\eqref{map} between forms on $M$ and semidensities on
$\Pi T^*M$ suggests that there exists a linear operator on
semidensities on odd symplectic supermanifolds corresponding to
the exterior differential.

\begin{df}[\cite{hov:semi}]
Let $\bs$ be a semidensity on an odd symplectic supermanifold $E$.
We assign to it a semidensity $\Delta {\bs}$ by the following
formula: if $\bs=s(x,\theta)\sqrt {{\mathcal D}(x,\theta)}$ in
some Darboux coordinates, then in the same coordinates
\begin{equation}\label{defoflaplaconsemidensity}
\Delta\bs:= \bigl(\Delta_0 s(x,\theta)\bigr)\sqrt{{\mathcal
D}(x,\theta)}= \frac{\p^2 s}{\p x^i\p\theta_i}\,\sqrt{\mathcal
D(x,\theta)}\,.
\end{equation}
We call this operator \textit{the canonical Laplacian on
semidensities}. (See \cite{hov:semi} for details.)
\end{df}

In the case of $E=\Pi T^*M$ one can see from~\eqref{map} that this
definition gives exactly the de Rham exterior differential:
\begin{equation}\label{commutation}
                  \omega\mapsto \bs_{\omega} \quad\Rightarrow\quad
                   \Delta\bs_{\omega}=\bs_{d{\omega}}\,.
\end{equation}
Of course, this relation is not a proof that  the operator given
by~\eqref{defoflaplaconsemidensity} is well-defined for a general
odd symplectic supermanifold $E$, because though one might
consider $E$ as $\Pi T^*M$ for some manifold $M$ the
identification~\eqref{map} fails under symplectomorphisms which
are not induced by diffeomorphisms of $M$.

To prove that the canonical operator is well-defined by
formula~\eqref{defoflaplaconsemidensity} one has to prove that the
r.h.s. of~\eqref{defoflaplaconsemidensity} indeed defines a
semidensity, i.e., under an arbitrary transformation from Darboux
coordinates $(x,\theta)$ to another Darboux coordinates
$(x',\theta')$ we have
\begin{equation}\label{correctness}
          \left(\Delta_0 s\right)\cdot
          \left(
          {\rm Ber}
            \frac{\p(x,\theta)}{\p(x',\theta')}
          \right)^{1/2}=
          \Delta_0^\prime \left(s\cdot\left({\rm Ber}
            \frac{\p(x,\theta)}{\p(x',\theta')}\right)^{1/2}
                   \right),
\end{equation}
where we denote by $\Delta_0$ and $\Delta^\prime_0$ the
``coordinate'' odd Laplacians~\eqref{defoflaplaciansimple} in the
Darboux coordinates $(x,\theta)$ and $(x',\theta')$, respectively.

Notice (see~\cite{hov:semi} for details) that every transformation
from Darboux coordinates $(x,\theta)$ to new Darboux coordinates
$(x',\theta')$ can be represented as the composition of
transformations of the following types:

\par\smallskip\noindent transformations corresponding
to $x'=x'(x)$:
\begin{equation}\label{diffcoord}
        x^\prime=x^\prime(x),\qquad \theta_{i'}=
                 \frac{\p x^i}{\p x^{i'}}\,\theta_i\,,
    \end{equation}
\par\smallskip\noindent   transformations  identical on the surface
$\theta=0$:
          \begin{equation}\label{adjustedcoord}
                   x'(x,\theta)\big\vert_{\theta=0}=x\,,\qquad
                   \theta'(x,\theta)\big\vert_{\theta=0}=\theta\,,
\end{equation}
\par\smallskip\noindent   transformations identical on even
coordinates:
\begin{equation}\label{specialcoord}
              x'=x,\qquad \theta_{i'}=\theta_i+\a_i\quad\text{such
              that}\ \
                \p_i\a_j-\p_j\a_i=0\,.
\end{equation}
It is sufficient to check  condition~\eqref{correctness} for
transformations~\eqref{diffcoord}, \eqref{adjustedcoord} and
\eqref{specialcoord} separately. For
transformations~\eqref{diffcoord} it follows from the
identity~\eqref{squarerelations}. The Berezinian of
transformation~\eqref{specialcoord} equals $1$,
hence~\eqref{correctness} is satisfied. One can show that
transformation~\eqref{adjustedcoord} is induced by a Hamiltonian
vector field (see \cite{hov:semi}), hence and it is sufficient to
check it infinitesimally.

Infinitesimal transformations are generated by  odd functions
(Hamiltonians) via the corresponding Hamiltonian vector fields. To
an odd Hamiltonian $Q(z)$ corresponds the infinitesimal canonical
transformation $\t z^A=z^A+\varepsilon \{Q,z^A\}$ generated by the
vector field ${\bf D}_Q$. The action of it on a semidensity
${\bs}$ can be expressed by a  ``differential''
$\delta_Q(s\sqrt{\D z})= \Delta_0 Q\cdot s\sqrt{\D z}-
\{Q,s\}\sqrt{\D z}$, because $\delta s=-\varepsilon\{Q,s\}$ and
$\delta \D z=\varepsilon\delta {\rm Ber}(\p z/\p\t z)\D
z=\varepsilon 2\Delta_0Q \D z$ for the infinitesimal
transformation generated by $Q$. Using $\Delta_0^2=0$ and
equation~\eqref{yieldingrelation}, we come to the commutation
relation $\Delta_0\delta_Q=\delta_Q\Delta_0$. Thus
condition~\eqref{correctness} is satisfied for infinitesimal
transformations.

\subsection {Properties of the canonical Laplacian}

The canonical Laplacian obviously obeys the condition
$\Delta^2=0$.

Let $\brho$ be an arbitrary volume form on an odd symplectic
supermanifold.  Then it is easy to check
using~\eqref{defoflaplaciana}, \eqref{yieldingrelation} and
\eqref{bracketpreservation} that the canonical Laplacian on
semidensities $\Delta$ obeys  the following condition:
\begin{equation}\label{newveryimportant}
    \Delta(f\sqrt\brho)=\left(\Delta_\brho
    f\right)\sqrt\brho+(-1)^{p(f)}f\Delta\sqrt\brho\,,
\end{equation}
where  $\Delta_\brho$ is the Laplacian~\eqref{defoflaplaciana} on
functions. Using \eqref{changingoflaplacian} one can rewrite this
relation in the following way, for the semidensity
$\bs=\sqrt\brho$:
\begin{equation}\label{voronovkommutator}
             [\Delta,f]\bs\equiv(\Delta\circ f-(-1)^{p(f)}f\circ\Delta)\bs=
             {\mathcal L}_{{\bf D}_f}\bs\,,
\end{equation}
where  ${\mathcal L}_{\bf D_f}$ is the Lie derivative along the
Hamiltonian vector field ${\bf D}_f$. Notice
that~\eqref{voronovkommutator} holds for an arbitrary semidensity
$\bs$, not only for an even non-degenerate semidensity
$\bs=\sqrt\brho$ corresponding to a volume form $\brho$. This
relation is very important for the study of Laplacians on
semidensities on arbitrary odd Poisson supermanifolds (see
\cite{sdens:voron}).

If $\brho$ is an arbitrary volume form, then by applying the
canonical Laplacian to the semidensity $\bs=\sqrt\brho$ we can
obtain a ``derived'' function $H=\Delta\sqrt{\brho}/\sqrt{\brho}$.
It turns out that the Hamiltonian vector field ${\bf D}_H$
corresponding to $H$ is nothing but the vector field
$\Delta_\brho^2$:
\begin{equation}\label{deltasquareb}
                    \Delta_\brho^2 f=
                    \left\{
                    \frac{\Delta\sqrt{\brho}}{\sqrt{\brho}}
                    ,f\right\}\,.
\end{equation}
(Compare this formula with~\eqref{cocyclecondition}.)

It is evident that if the  form $\brho$ is normal, then by
definition
\begin{equation}\label{batalinidentityaplication}
                                   \Delta\sqrt\brho=0\,.
\end{equation}
This is just an invariant expression for the Batalin--Vilkovisky
identity~\eqref{batalinindentity}

Now let us return to the relation between differential forms on
$M$ and semidensities on $\Pi T^*M$, and to its  generalization
for arbitrary odd symplectic supermanifolds.

We call a semidensity $\bs$ \textit{closed} or \textit{exact} if
$\Delta\bs=0$ or $\bs=\Delta{{\boldsymbol r}}$ respectively. The
condition $\Delta^2=0$ for the canonical Laplacian corresponds to
$d^2=0$ for the exterior differential.

Equations~\eqref{map}, \eqref{commutation} and
\eqref{voronovkommutator} allow the translation of formulae of
vector calculus on $M$ into  formulae for semidensities on $\Pi
T^*M$ (see \cite{hov:semi}). For example, under the map
\eqref{map} the ``interior multiplication'' of a differential form
$\omega$ by a multivector field $T$ transforms into the usual
product of the semidensity ${\bs}_{\omega}$ with the function
corresponding to the multivector field. Hence,
equation~\eqref{voronovkommutator} corresponds to the formula for
the Lie derivative of a differential form along a multivector
field (a generalization of Cartan's homotopy formula).

Consider in more details  the following two constructions which do
not appear naturally in classical calculus of forms. (Below we use
familiar formal properties of the Fourier transform.)

a)  If $a=a_i(x)dx^i$  is a $1$-form on $M$ and a semidensity
${\bs}=s(x,\theta)\sqrt{{\mathcal D}(x,\theta)}$ on $\Pi T^*M$
corresponds to another form $\omega$,  then one can see that the
semidensity $a_i{\p s\over\p\theta_i}\sqrt{D(x,\theta)}$
corresponds to the form $a\wedge \omega$. (Compare with the
familiar relation between the differentiation and multiplication
by a coordinate for the classical Fourier transform.) Consider the
following generalization. Let $a=a_idx^i$ be a one-form on $M$
with odd coefficients (we have to allow ``external odd
parameters'' for this). For an arbitrary semidensity
${\bs}=s(x,\theta)\sqrt{D(x,\theta)}$ consider a new semidensity
${\bs}^\prime$, which we denote by $a\,{\mathcal d}\,{\bs}$, given
by the formula ${\bs}^\prime= a\,{\mathcal
d}\,{\bs}:=s(x,\theta_i+a_i)\sqrt{D(x,\theta)}$. It is a
well-defined operation, because the coefficients $a_i$ have the
same transformation law as the variables $\theta_i$. (Notice that
the Berezinian of the transformation $(x^i,\theta_i)\mapsto
(x^i,\theta_i+a_i)$ equals $1$.)   Respectively, if the
semidensity ${\bs}$ corresponds to a differential form
$\omega=\sum \omega_k$, then we denote by $a\,{\mathcal
d}\,\omega$ the differential form such that the semidensity
$a\,{\mathcal d}\,{\bs}$ corresponds to $a\,{\mathcal d}\,\omega$.
One can see that
\begin{equation*}
    a\,{\mathcal d}\,\omega=\sum_{p=0}^k {1\over p!}
    \underbrace {a\wedge\dots\wedge a}_{p\quad{\rm times}}\wedge
                          \omega_{k-p}\,,\quad (k=0,\dots,n)\,.
\end{equation*}
We obtain an action of the abelian supergroup of differential
one-forms ``with odd values'' (i.e., $\Pi \Omega^1(M)$) in the
spaces of semidensities and differential forms (see
\cite{hov:semi}).

b) Let $\omega=\sum \omega_k$ and $\omega^\prime=\sum
\omega^\prime_k$ be differential forms on $M^n$ such that their
top-degree components $\omega_n$ and $\omega^\prime_n$ are
non-zero, and let the semidensities ${\bs}$ and ${\bs^\prime}$
correspond to $\omega$ and $\omega^\prime$ respectively. Then we
can define a new form $\tilde \omega:=\omega*\omega^\prime$ such
that the corresponding semidensity is equal to $\sqrt{{\bs}\cdot
{\bs}^\prime}$. The condition $\omega_n\not=0$,
$\omega^\prime_n\not=0$ for the top-degree components makes the
square root operation uniquely defined.

The 1-1 correspondence between forms on $M$ and semidensities on
$\Pi T^*M$ is defined using  the  cotangent bundle structure on
the odd symplectic supermanifold $\Pi T^*M$. Bearing in mind that
every odd symplectic manifold $E$  with an underlying manifold $M$
is symplectomorphic to $\Pi T^*M$ (see \cite{hov:semi} and the
footnote in subsection~\ref{subsecdefinitionandproperties}), let
us analyze the  relation between semidensities on $E$ and
differential forms on $M$. The map~\eqref{map} is not invariant
under arbitrary symplectomorphisms of the total symplectic
supermanifold $E=\PTM$. In other words, if $L$ is an arbitrary
$(n.0)$-dimensional Lagrangian surface in $E$, then  the
correspondence between semidensities on $E$ and differential forms
on $L$ depends on an \textit{identifying symplectomorphism}, i.e.
a symplectomorphism $\varphi:\, \Pi T^*L\to E$ such that
$\varphi\vert_L={\mathrm{id}}$.

Consider the  following symplectomorphisms of an odd symplectic
supermanifold $E=\PTM$:
\begin{equation}\label{diff}
          \text{symplectomorphisms induced by diffeomorphisms of $M$}
\end{equation}
(these symplectomorphisms preserve the cotangent bundle
structure), symplectomorphisms ``adjusted to $M$'', i.e. identical
on $M$:
\begin{equation}\label{adjusted}
   \varphi\colon\quad \varphi^*\omega=\omega\,, \ \varphi\big\vert_{M}=\text{$\mathrm{id}$}
\end{equation}
(they destroy the cotangent bundle structure except for
$\varphi=\mathrm{id}$), symplectomorphisms corresponding to closed
one-forms on $M$:
\begin{equation}\label{special}
           \varphi^*x^{i}=x^i\,,\
           \varphi^*\theta_{i}=\theta_i+\alpha_i(x)\,,\
                     (\p_i\a_j-\p_j\a_i=0)
\end{equation}
where $(x^i,\theta_i)$ are coordinates on $\Pi T^*M$ corresponding
to some coordinates $(x^i)$  on $M$ and $\alpha=\alpha_i(x)dx^i$
is a closed one-form on $M$ with odd values (these
symplectomorphisms move the Lagrangian surface $M$; notice that an
arbitrary $(n.0)$-dimensional Lagrangian surface $L$ is given by
the equations $\theta_i-\a_i(x)=0$ where $\a_i(x)dx^i$ is a closed
odd-valued one-form). It might be worth noting that
symplectomorphisms of $\PTM$ form a supergroup.

One can prove that an arbitrary symplectomorphism of $\PTM$ can be
represented as the composition of symplectomorphisms \eqref{diff},
\eqref{adjusted} and \eqref{special} (see \cite{hov:semi}).
(Compare this statement with the statement that an arbitrary
transformation of Darboux coordinates can be represented as the
composition of transformations \eqref{diffcoord},
\eqref{adjustedcoord}  and \eqref{specialcoord}).

Notice that every adjusted symplectomorphism \eqref{adjusted} has
the following appearance:
\begin{equation}\label{appearance}
                           \begin{gathered}
                            x^i\rightarrow x^i+f^i(x,\theta)\,,\quad
                                {\rm where}\,\, f^i(x,\theta)=O(\theta)\\
                                \theta_i\rightarrow \theta_i+g_i(x,\theta)\,,\quad
                                {\rm where}\,\, g_i(x,\theta)=O(\theta^2)\\
                           \end{gathered}
\end{equation}
where $(x^i,\theta_i)$ are the Darboux coordinates on $\PTM$
corresponding to coordinates $x^i$ on $M$. One can show that there
exists a Hamiltonian
$Q(x,\theta)=Q^{ik}(x,\theta)\theta_i\theta_k$ that generates this
transformation, i.e.,  \eqref{appearance} can be included in a
$1$-parameter family of transformations defined by the
differential equation $\dot z=\{Q,z\}$ (see \cite{hov:semi} for
details).

Comparing~\eqref{appearance} with~\eqref{map} we arrive at an
important conclusion:

\begin{prop}\label{proptopcomp}
The top-degree component of the form corresponding to a
semidensity on  $\PTM$ does not change under any symplectomorphism
adjusted to $M$. In other words, a  semidensity on an odd
symplectic supermanifold $E$ defines a volume form (density) for
all $(n.0)$-dimensional Lagrangian surfaces\,\footnote {A relation
between semidensities on $E$ and densities on Lagrangian surfaces
can be defined for arbitrary Lagrangian surfaces (see
\cite{ass:bv} and \cite{hov:khn1}).}.
\end{prop}

Now consider a \textit{closed} semidensity on $\PTM$. To it
corresponds a closed differential form on $M$. It follows
from~\eqref{voronovkommutator} and~\eqref{adjusted} that the
action  of an adjusted symplectomorphism \eqref{appearance}
changes the semidensity and the corresponding form by an exact
semidensity and an exact form respectively. Hence, we arrive at
another important conclusion:

\begin{prop}\label{propcohclass}
To a closed semidensity on $\PTM$ corresponds a cohomology class
of differential forms on $M$ independently of the bundle
structure. If two closed semidensities $\bs$ and $\bs^\prime$
coincide on $M$ and differ by an exact semidensity, then there
exists an adjusted symplectomorphism $\varphi:\, \PTM\to \PTM$
such that $\varphi^*\bs=\bs^\prime$.
\end{prop}

Propositions~\ref{proptopcomp} and ~\ref{propcohclass} were stated
and proved in  \cite{ass:bv} and  in \cite{hov:semi}, but in the
work~\cite{ass:bv}  semidensities do not appear explicitly.

Based on the concept of semidensities, the properties of the
canonical Laplacian and the above Propositions we will now analyze
the  statements \eqref{nashiteoremy} concerning the
Batalin--Vilkovisky formalism.

\subsection {Master  equation on semidensities}

The claim is that \textit{the Batalin--Vilkovisky master equation
\eqref{masterequation} is an equation on the semidensity
$\bs=\sqrt {f(x,\theta)\brho}$. A solution of the
Batalin--Vilkovisky quantum master equation is a closed
semidensity: $\,\Delta\bs=0$.}

Suppose $E$ is an odd symplectic supermanifold with the compact
connected orientable underlying manifold $M$. $E$ can be
identified with $\Pi T^*L$ for every closed $(n.0)$-dimensional
Lagrangian submanifold $L$ (see~\cite{hov:semi} for details) and
any two identifications differ by an adjusted symplectomorphism.
An arbitrary $(n.0)$-dimensional closed Lagrangian surface $L$ is
given by a closed one-form on $M$ (see \eqref{special}).

Let us now rewrite the implications~\eqref{nashiteoremy} in terms
of semidensities:
\begin{equation}\label{nashiteoremyluchshe}
    \text{$\brho$ is a normal volume form}\ \ \Rightarrow\ \
       \Delta \sqrt{\brho}=0\ \ \Rightarrow\ \
                  \Delta_\brho^2=0\,.
\end{equation}
The first implication follows  from the definition of the
canonical operator $\Delta$. The second implication  follows from
equation~\eqref{deltasquareb}. Let us analyze to what extent these
conditions are equivalent. Let $\brho$ be a volume form such that
$\Delta_\rho^2=0$. By~\eqref{deltasquareb}, then it follows that
$\Delta\sqrt\brho= \nu \sqrt\rho$, where $\nu$ is an odd constant.
(If external odd parameters are not allowed, then, of course,
$\nu=0$. Our analysis takes into consideration possible ``odd
moduli''.) This odd constant is the obstruction to the condition
$\Delta\sqrt\brho=0$, i.e. to the closedness of the semidensity
$\sqrt\brho$.

Suppose $\nu=0$.  Then the master equation $\Delta\sqrt\brho=0$ is
satisfied. Consider in this case an arbitrary closed
$(n.0)$-dimensional Lagrangian surface $L$ and an arbitrary
identifying symplectomorphism $\varphi\colon\, E\to \Pi T^*L$.
Under the map~\eqref{map} to every closed semidensity $\sqrt\rho$
on $E$ corresponds a closed differential form
$\omega=\omega_0+\omega_1+\dots+\omega_n$ on $L$, where the top
degree form $\omega_n$ defines a volume form on $L$. The closed
$0$-form $\omega_0$ is a constant. It is easy to see
from~\eqref{map} and \eqref{special} that the value of this
constant (up to a sign) does not depend on the choice of the
Lagrangian surface and on the choice of the identifying
symplectomorphism. The top degree form, clearly,  depends on the
Lagrangian surface but does not depend on the identifying
symplectomorphism (by Proposition~\ref{proptopcomp}).  On the
other hand, all other closed forms $\omega_k$ ($1<k<n$) can be
eliminated by a suitable choice of the identifying
symplectomorphism.

We have arrived, finally, to the following theorem:

\begin{thm}
Let $E^{n.n}$ be an odd symplectic supermanifold with the closed
orientable compact underlying manifold $M$. Let $\rho$ be a volume
form on $E^{n.n}$ such that $\Delta_\rho^2=0$.  To this volume
form corresponds an odd constant $\nu$. If this odd constant is
equal to zero, then the volume form defines a closed semidensity
$\bs=\sqrt\rho$, a solution of the Batalin--Vilkovisky quantum
master equation. To a closed semidensity corresponds a constant
$c$ defined by the zero cohomology class of the differential form
corresponding to the semidensity $\bs$. If this constant is equal
to zero, then the volume form $\brho$ is normal.
\end{thm}

\section{Discussion}

The existence of Darboux coordinates that can locally make flat
every surface in an odd symplectic supermanifold together with the
absence of an invariant volume form   make  odd symplectic
geometry a poor candidate for finding local invariants if no extra
structure is provided. Hence the existence of the canonical odd
Laplacian~\eqref{defoflaplaconsemidensity}  looks mysterious. We
shall to explain this fact briefly . (See details
in~\cite{sdens:voron})

Consider an arbitrary $n$-th order linear operator acting on
functions or densities of some weight $t$ on an arbitrary manifold
(or supermanifold). One can consider its principal symbol, i.e.,
the coefficients at the highest order derivatives. It is a
contravariant tensor field of rank $n$.  In the case of an $n$-th
order operator $\hat A$ acting on  densities of  weight $t$, the
adjoint operator ${\hat{A^\dagger}}$ acting on  densities of
weight $1-t$ can be defined by the equation $\int\bs_1\cdot(\hat
A\bs_2)=\int({\hat{A^\dagger}}\bs_1)\cdot\bs_2$, where $\bs_2$ is
an arbitrary density of  weight $t$ and $\bs_1$ is an arbitrary
density of  weight $1-t$.  Hence in the case of an operator $\hat
A$ acting on semidensities ($t=1/2$) the operators $\hat A$ and
$\hat A^\dagger$ act on the same space. Assuming that the
coefficients are real, the operators $\hat A$ and $\hat A^\dagger$
have the same principal symbol. One can consider the principal
symbol of $\hat A-\hat A^\dagger(-1)^n$, which is a tensor field
of  rank $n-1$.  It is the so-called \textit{subprincipal symbol}
of the operator $\hat A$.

Let us consider the canonical operator $\Delta$ on an odd
symplectic supermanifold in arbitrary coordinates (not necessary
Darboux coordinates). The highest order coefficients make  the
principal symbol, which is here the tensor of  rank $2$ defining
the odd symplectic structure.  (More precisely it is the tensor
$\Sch^{AB}$ that defines the master
Hamiltonian~\eqref{masterodd}.) It is easy to see that if $\hat A$
is an arbitrary linear differential operator  of the second order
on semidensities having the principal symbol defined by the odd
symplectic structure,  then the subprincipal symbol of this
operator is equal to $[\hat A,f]-{\mathcal L}_{{\bf D}_f}$. Hence
it follows from \eqref{voronovkommutator} that the subprincipal
symbol of the canonical Laplacian $\Delta$ is equal to zero. The
coefficients at the first derivatives are fixed by this condition.
In fact, these two conditions on the principal symbol and
subprincipal symbol are equivalent to the
equation~\eqref{voronovkommutator}. An arbitrary linear operator
$\Delta^\prime$ on semidensities obeying   the
condition~\eqref{voronovkommutator} is equal to $\Delta+C$, where
$C$ is a scalar (a zero-order operator). The condition
$\Delta^\prime\sqrt\brho=0$ for an arbitrary normal volume form
fixes this scalar $C=0$ according
to~\eqref{batalinidentityaplication}: hence, \textit{on an odd
symplectic supermanifold there is no distinguished volume form,
but there is a distinguished class of normal volume forms.}

\medskip \textbf{Acknowledgment:} I am deeply grateful to
Th.~Voronov for encouraging me to write this paper and for
essential editorial help.

\end{document}